\documentclass[12pt]{article}

\catcode`\@=11

\thicklines
\newskip\Einheit \Einheit=.6cm
\newcount\xcoord \newcount\ycoord
\newdimen\xdim \newdimen\ydim \newdimen\PfadD@cke \newdimen\Pfadd@cke
\def\PfadDicke#1{\PfadD@cke#1 \divide\PfadD@cke by2 
\Pfadd@cke\PfadD@cke \multiply\PfadD@cke by2}
\long\def\LOOP#1\REPEAT{\def\BODY{#1}\ITERATE}
\def\ITERATE{\BODY \let\next\ITERATE \else\let\next\relax\fi \next}
\let\REPEAT=\fi
\def\Punkt{\hbox{\raise-2pt\hbox to0pt{\hss\scriptsize$\bullet$\hss}}}

\def\DuennPunkt(#1,#2){\unskip
  \raise#2 \Einheit\hbox to0pt{\hskip#1 \Einheit
          \raise-1.5pt\hbox to0pt{\hss\tiny$\bullet$\hss}\hss}}
		  
\def\NormalPunkt(#1,#2){\unskip
  \raise#2 \Einheit\hbox to0pt{\hskip#1 \Einheit
          \raise-3pt\hbox to0pt{\hss\large$\bullet$\hss}\hss}}
\def\DickPunkt(#1,#2){\unskip
  \raise#2 \Einheit\hbox to0pt{\hskip#1 \Einheit
          \raise-4pt\hbox to0pt{\hss\Large$\bullet$\hss}\hss}}
\def\Kreis(#1,#2){\unskip
  \raise#2 \Einheit\hbox to0pt{\hskip#1 \Einheit
          \raise-4pt\hbox to0pt{\hss\Large$\circ$\hss}\hss}}
\def\Diagonale(#1,#2)#3{\unskip\leavevmode
  \xcoord#1\relax \ycoord#2\relax
      \raise\ycoord \Einheit\hbox to0pt{\hskip\xcoord \Einheit
         \unitlength\Einheit
         \line(1,1){#3}\hss}}
\def\AntiDiagonale(#1,#2)#3{\unskip\leavevmode
  \xcoord#1\relax \ycoord#2\relax \advance\xcoord by -0.05\relax
      \raise\ycoord \Einheit\hbox to0pt{\hskip\xcoord \Einheit
         \unitlength\Einheit
         \line(1,-1){#3}\hss}}
\def\Pfad(#1,#2),#3\endPfad{\unskip\leavevmode
  \xcoord#1 \ycoord#2 \thicklines\ZeichnePfad#3\endPfad\thinlines}
\def\ZeichnePfad#1{\ifx#1\endPfad\let\next\relax
  \else\let\next\ZeichnePfad
    \ifnum#1=1
      \raise\ycoord \Einheit\hbox to0pt{\hskip\xcoord \Einheit
         \vrule height\Pfadd@cke width1 \Einheit depth\Pfadd@cke\hss}%
      \advance\xcoord by 1
     \else\ifnum#1=2
      \raise\ycoord \Einheit\hbox to0pt{\hskip\xcoord \Einheit
         \unitlength\Einheit
         \line(0,1){1}\hss}
      \advance\xcoord by 0
      \advance\ycoord by 1
 \else\ifnum#1=3
      \raise\ycoord \Einheit\hbox to0pt{\hskip\xcoord \Einheit
         \unitlength\Einheit
         \line(1,1){1}\hss}
      \advance\xcoord by 1
      \advance\ycoord by 1
    \else\ifnum#1=4
      \raise\ycoord \Einheit\hbox to0pt{\hskip\xcoord \Einheit
         \unitlength\Einheit
         \line(1,-1){1}\hss}
      \advance\xcoord by 1
      \advance\ycoord by -1
   \else\ifnum#1=5
      \raise\ycoord \Einheit\hbox to0pt{\hskip\xcoord \Einheit
         \unitlength\Einheit
         \line(2,1){2}\hss}
      \advance\xcoord by 2
      \advance\ycoord by 1
	  \else\ifnum#1=6
      \raise\ycoord \Einheit\hbox to0pt{\hskip\xcoord \Einheit
         \unitlength\Einheit
         \line(2,-1){2}\hss}
      \advance\xcoord by 2
      \advance\ycoord by -1
	  \else\ifnum#1=7
      \raise\ycoord \Einheit\hbox to0pt{\hskip\xcoord \Einheit
         \unitlength\Einheit
         \line(3,1){3}\hss}
      \advance\xcoord by 3
      \advance\ycoord by 1
	  \else\ifnum#1=8
      \raise\ycoord \Einheit\hbox to0pt{\hskip\xcoord \Einheit
         \unitlength\Einheit
         \line(3,-1){3}\hss}
      \advance\xcoord by 3
      \advance\ycoord by -1
    \fi\fi\fi\fi\fi\fi\fi\fi
  \fi\next}
\def\hSSchritt{\leavevmode\raise-.4pt\hbox 
to0pt{\hss.\hss}\hskip.2\Einheit
  \raise-.4pt\hbox to0pt{\hss.\hss}\hskip.2\Einheit
  \raise-.4pt\hbox to0pt{\hss.\hss}\hskip.2\Einheit
  \raise-.4pt\hbox to0pt{\hss.\hss}\hskip.2\Einheit
  \raise-.4pt\hbox to0pt{\hss.\hss}\hskip.2\Einheit}
\def\vSSchritt{\vbox{\baselineskip.2\Einheit\lineskiplimit0pt
\hbox{.}\hbox{.}\hbox{.}\hbox{.}\hbox{.}}}
\def\DSSchritt{\leavevmode\raise-.4pt\hbox to0pt{%
  \hbox to0pt{\hss.\hss}\hskip.2\Einheit
  \raise.2\Einheit\hbox to0pt{\hss.\hss}\hskip.2\Einheit
  \raise.4\Einheit\hbox to0pt{\hss.\hss}\hskip.2\Einheit
  \raise.6\Einheit\hbox to0pt{\hss.\hss}\hskip.2\Einheit
  \raise.8\Einheit\hbox to0pt{\hss.\hss}\hss}}
\def\dSSchritt{\leavevmode\raise-.4pt\hbox to0pt{%
  \hbox to0pt{\hss.\hss}\hskip.2\Einheit
  \raise-.2\Einheit\hbox to0pt{\hss.\hss}\hskip.2\Einheit
  \raise-.4\Einheit\hbox to0pt{\hss.\hss}\hskip.2\Einheit
  \raise-.6\Einheit\hbox to0pt{\hss.\hss}\hskip.2\Einheit
  \raise-.8\Einheit\hbox to0pt{\hss.\hss}\hss}}
\def\SPfad(#1,#2),#3\endSPfad{\unskip\leavevmode
  \xcoord#1 \ycoord#2 \ZeichneSPfad#3\endSPfad}
\def\ZeichneSPfad#1{\ifx#1\endSPfad\let\next\relax
  \else\let\next\ZeichneSPfad
    \ifnum#1=1
      \raise\ycoord \Einheit\hbox to0pt{\hskip\xcoord \Einheit
         \hSSchritt\hss}%
      \advance\xcoord by 1
    \else\ifnum#1=2
      \raise\ycoord \Einheit\hbox to0pt{\hskip\xcoord \Einheit
        \hbox{\hskip-2pt \vSSchritt}\hss}%
      \advance\ycoord by 1
    \else\ifnum#1=3
      \raise\ycoord \Einheit\hbox to0pt{\hskip\xcoord \Einheit
         \DSSchritt\hss}
      \advance\xcoord by 1
      \advance\ycoord by 1
    \else\ifnum#1=4
      \raise\ycoord \Einheit\hbox to0pt{\hskip\xcoord \Einheit
         \dSSchritt\hss}
      \advance\xcoord by 1
      \advance\ycoord by -1
    \fi\fi\fi\fi
  \fi\next}
\def\Koordinatenachsen(#1,#2){\unskip
 \hbox to0pt{\hskip-.5pt\vrule height#2 \Einheit width.5pt depth1 
\Einheit}%
 \hbox to0pt{\hskip-1 \Einheit \xcoord#1 \advance\xcoord by1
    \vrule height0.25pt width\xcoord \Einheit depth0.25pt\hss}}
\def\Koordinatenachsen(#1,#2)(#3,#4){\unskip
 \hbox to0pt{\hskip-.5pt \ycoord-#4 \advance\ycoord by1
    \vrule height#2 \Einheit width.5pt depth\ycoord \Einheit}%
 \hbox to0pt{\hskip-1 \Einheit \hskip#3\Einheit 
    \xcoord#1 \advance\xcoord by1 \advance\xcoord by-#3 
    \vrule height0.25pt width\xcoord \Einheit depth0.25pt\hss}}
\def\Gitter(#1,#2){\unskip \xcoord0 \ycoord0 \leavevmode
  \LOOP\ifnum\ycoord<#2
    \loop\ifnum\xcoord<#1
      \raise\ycoord \Einheit\hbox to0pt{\hskip\xcoord 
\Einheit\Punkt\hss}%
      \advance\xcoord by1
    \repeat
    \xcoord0
    \advance\ycoord by1
  \REPEAT}
\def\Gitter(#1,#2)(#3,#4){\unskip \xcoord#3 \ycoord#4 \leavevmode
  \LOOP\ifnum\ycoord<#2
    \loop\ifnum\xcoord<#1
      \raise\ycoord \Einheit\hbox to0pt{\hskip\xcoord 
\Einheit\Punkt\hss}%
      \advance\xcoord by1
    \repeat
    \xcoord#3
    \advance\ycoord by1
  \REPEAT}
\def\Label#1#2(#3,#4){\unskip \xdim#3 \Einheit \ydim#4 \Einheit
  \def\lo{\advance\xdim by-.5 \Einheit \advance\ydim by.5 \Einheit}%
  \def\llo{\advance\xdim by-.25cm \advance\ydim by.5 \Einheit}%
  \def\loo{\advance\xdim by-.5 \Einheit \advance\ydim by.25cm}%
  \def\o{\advance\ydim by.25cm}%
  \def\ro{\advance\xdim by.5 \Einheit \advance\ydim by.5 \Einheit}%
  \def\rro{\advance\xdim by.25cm \advance\ydim by.5 \Einheit}%
  \def\roo{\advance\xdim by.5 \Einheit \advance\ydim by.25cm}%
  \def\l{\advance\xdim by-.30cm}%
  \def\r{\advance\xdim by.30cm}%
  \def\lu{\advance\xdim by-.5 \Einheit \advance\ydim by-.6 \Einheit}%
  \def\llu{\advance\xdim by-.25cm \advance\ydim by-.6 \Einheit}%
  \def\luu{\advance\xdim by-.5 \Einheit \advance\ydim by-.30cm}%
  \def\u{\advance\ydim by-.30cm}%
  \def\ru{\advance\xdim by.5 \Einheit \advance\ydim by-.6 \Einheit}%
  \def\rru{\advance\xdim by.25cm \advance\ydim by-.6 \Einheit}%
  \def\ruu{\advance\xdim by.5 \Einheit \advance\ydim by-.30cm}%
  #1\raise\ydim\hbox to0pt{\hskip\xdim
     \vbox to0pt{\vss\hbox to0pt{\hss$#2$\hss}\vss}\hss}%
}
\catcode`\@=12

\usepackage{amsmath,amsthm,amscd}
\usepackage{pstricks,pst-node}
\usepackage{amssymb}
\usepackage{color}
\usepackage{xspace}
\usepackage[colorlinks=true,
linkcolor=green,
filecolor=brown,
citecolor=green]{hyperref}

\def\red{\textcolor{red} }
\def\blue{\textcolor{blue} }
\def\v{\vert}
\def\ep{\epsilon}
\def\gl{ground level}
\def\gf{generating function}
\def\a{\ensuremath{\mathbf p}\xspace}
\def\c{\ensuremath{\mathbf c}\xspace}
\def\d{\ensuremath{\mathbf d}\xspace}
\def\e{\ensuremath{\mathbf e}\xspace}
\def\vv{\ensuremath{\mathbf v}\xspace}
\def\u{\ensuremath{\mathbf p}\xspace}
\def\cc{\ensuremath{\mathcal C}\xspace}
\def\p{\ensuremath{\mathcal P}\xspace}
\def\ss{\ensuremath{S}\xspace}
\def\s{\ensuremath{\mathbf S}\xspace}

\begin{document}

\newtheorem{theorem}{Theorem}
\newtheorem*{defn}{Definition}
\newtheorem{lemma}[theorem]{Lemma}
\newtheorem{prop}[theorem]{Proposition}
\newtheorem{cor}[theorem]{Corollary}
\begin{center}
{\Large
A Bijection on Dyck Paths and Its Cycle Structure
}

\vspace{10mm}
DAVID CALLAN  \\
Department of Statistics  \\
\vspace*{-1mm}
University of Wisconsin-Madison  \\
\vspace*{-1mm}
1300 University Ave  \\
\vspace*{-1mm}
Madison, WI \ 53706-1532  \\
{\bf callan@stat.wisc.edu}  \\
\vspace{5mm}

November 21, 2006 
\end{center}

\vspace{3mm}
\begin{abstract}
The known bijections on Dyck paths are either involutions or have 
notoriously intractable cycle structure. Here we present a size-preserving 
bijection on Dyck paths whose cycle structure is amenable to complete analysis.
In particular, each cycle has length a power of 2. 
A new manifestation of the 
Catalan numbers as labeled forests crops up enroute as does the Pascal matrix mod 2. We use the bijection 
to show the equivalence of two known manifestations 
of the Motzkin numbers.
Finally, we consider some statistics on the new Catalan manifestation. 
\end{abstract}

\vspace{10mm}

{\Large \textbf{1 \ Introduction}  }\quad
There are several bijections on Dyck paths in the literature \cite{twobij04,catfine,
invol1999,bij1998,ordered,simple2003,don80,acp44,lalanne92,lalanne93,vaille97}, usually 
introduced to show the equidistribution of statistics: if a bijection 
sends statistic A to statistic B, then clearly both have the same 
distribution. Another aspect of such a bijection is its cycle structure considered as 
a permutation on Dyck paths. Apart from involutions, this question is usually 
intractable. For example, Donaghey \cite{don80} introduces a 
bijection, gets some results on a restriction version, and notes its apparently chaotic behavior in general.  
In similar vein, Knuth \cite{acp44} defines a conjugate ($R$)  
and transpose ($T$), both involutions, on ordered forests, equivalently on 
Dyck paths, and asks when they commute \cite[Ex.\,17,\,7.2.1.6]{acp44}, equivalently, what are the 
fixed points of $(RT)^{2}$? This question is still open. (Donaghey's 
bijection is equivalent to the composition $RT$.)

In this paper, after reviewing Dyck path terminology (\S2), we 
recursively define a new bijection $F$ on Dyck paths (\S 3) and analyze its cycle structure (\S4,\:\S5). 
\S 4 treats the restriction of $F$ to paths that avoid 
the subpath $DUU$, and involves an encounter with the Pascal matrix mod 2.
\S 5 generalizes to arbitrary paths. This entails an explicit description of $F$ involving
a new manifestation of the Catalan numbers as certain colored 
forests in which each vertex is labeled with an integer composition.
We show that each orbit has length a power of 2, find
generating functions for orbit size, and characterize paths with given orbit 
size in terms of subpath avoidance. In particular, the 
fixed points of $F$ are those Dyck paths that avoid $DUDD$ 
and $UUP^{+}DD$ where $P^{+}$ denotes a 
nonempty Dyck path. 
\S 6 uses the bijection $F$ to show the equivalence of two known manifestations 
of the Motzkin numbers.
\S7 considers some statistics on the new Catalan manifestation.

\vspace{10mm}

{\Large \textbf{2 \ Dyck Path Terminology}  }\quad
A Dyck path, as usual, is a lattice path of upsteps $U=(1,1)$ and 
downsteps $D=(1,-1)$, the same number of each, that stays weakly 
above the horizontal line joining its initial and terminal points (vertices). A peak is an occurrence of $UD$, a valley is a 
$DU$. 

\vspace*{-3mm}

\Einheit=0.6cm
\[
\Label\o{\rightarrow}(-5.2,3)
\Label\u{\uparrow}(-2,1.9)
\Label\l{ \textrm{{\footnotesize peak upstep}}}(-6.9,3.5)
\Label\u{ \textrm{{\footnotesize valley}}}(-2,1)
\Label\u{ \textrm{{\footnotesize vertex}}}(-2,0.4)
\Label\u{\uparrow}(.4,1.5)
\Label\u{ \textrm{{\footnotesize return}}}(.4,.8)
\Label\u{ \textrm{{\footnotesize downstep}}}(.4,0.2)
\SPfad(-7,1),1111\endSPfad
\SPfad(1,1),111111\endSPfad
\SPfad(-1,1),1\endSPfad
\Pfad(-7,1),33344344334344\endPfad
\DuennPunkt(-7,1)
\DuennPunkt(-6,2)
\DuennPunkt(-5,3)
\DuennPunkt(-4,4)
\DuennPunkt(-3,3)
\DuennPunkt(-2,2)
\DuennPunkt(-1,3) 
\DuennPunkt(0,2) 
\DuennPunkt(1,1)
\DuennPunkt(2,2)
\DuennPunkt(3,3)
\DuennPunkt(4,2)
\DuennPunkt(5,3)
\DuennPunkt(6,2)
\DuennPunkt(7,1)
\Label\u{\uparrow}(4,1)
\Label\u{ \textrm{{\footnotesize \gl}}}(4,.2)
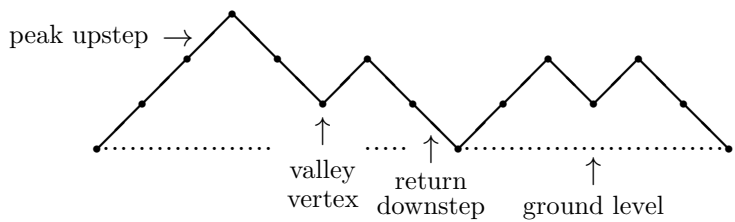
\Label\o{ \textrm{\small  A Dyck 7-path with 2 components, 2$DUD$s, and height 3}}(0,-2.5)
\]

\vspace*{1mm}

The size (or semilength) of a Dyck path is its number of upsteps and a 
Dyck path of size $n$ is a Dyck $n$-path. The empty Dyck path (of size 
0) is denoted $\ep$. The number of Dyck $n$-paths is the Catalan 
number $C_{n}$, sequence 
\htmladdnormallink{A000108}{http://www.research.att.com:80/cgi-bin/access.cgi/as/njas/sequences/eisA.cgi?Anum=A000108}
in 
\htmladdnormallink{OEIS}{http://www.research.att.com/~njas/sequences/Seis.html} .
The height of a vertex in a Dyck path is its vertical height above 
\gl\ and the height of the path is the maximum height of its vertices. 
A return downstep is one that returns the path to \gl. A \emph{primitive} 
Dyck path is one with exactly one return (necessarily at the end). 
Note that the empty Dyck path $\ep$ is not primitive. Its returns split 
a nonempty Dyck path into one or more primitive Dyck paths, called its 
\emph{components}. Upsteps and downsteps come in matching pairs: travel due east 
from an upstep to the first downstep encountered. More precisely, 
$D_{0}$ is the matching downstep for upstep $U_{0}$ if $D_{0}$ 
terminates the shortest Dyck subpath that starts with $U_{0}$. 
We use \p to denote the set of primitive Dyck paths, $\p_{n}$ for 
$n$-paths, $\p(DUU)$ for those that avoid $DUU$ as a subpath, and 
$\p[DUU]$ for those that contain at least one $DUU$. A path 
$UUUDUDDD$, for example, is abbreviated $U^{3}DUD^{3}$.

\vspace{10mm}

{\Large \textbf{3 \ The Bijection}  }\quad 
Define a size-preserving bijection $F$ on Dyck paths recursively as follows.
First, $F(\ep)=\ep$ and for a non-primitive Dyck path $P$ with components 
$P_{1},P_{2},\ldots,P_{r}\ (r\ge 2)$, $F(P)=F(P_{1})F(P_{2})\ldots 
F(P_{r})$ (concatenation). This reduces matters to primitive paths. 
From a consideration of the last vertex at height 3 (if any), every primitive Dyck path $P$ has the 
form $UQ(UD)^{i}D$ with $i\ge 0$ and $Q$ a Dyck path that 
is either empty (in case no vertex is at height 3) or ends $DD$; define $F(P)$ by
\[
F(P)=
\begin{cases}
    U^{i}F(R)UDD^{i} & \textrm{if $Q$ is primitive, say $Q=URD$, and} \\
    U^{i+1}F(Q)D^{i+1} & \textrm{if $Q$ is not primitive.}
\end{cases}
\]
Schematically,

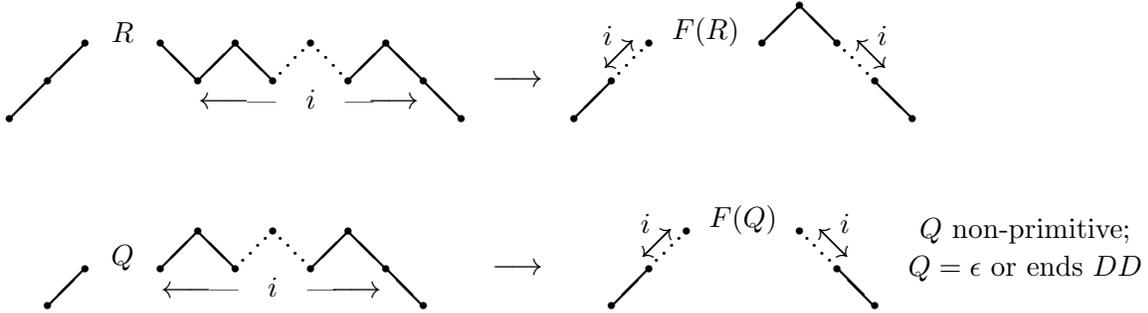
\vspace*{-5mm}

\Einheit=0.5cm
\[
\Pfad(-15,5),33\endPfad
\Pfad(-11,7),434\endPfad
\Pfad(-6,6),344\endPfad
\Pfad(0,5),3\endPfad
\Pfad(5,7),34\endPfad
\Pfad(8,6),4\endPfad
\SPfad(-8,6),34\endSPfad
\SPfad(1,6),3\endSPfad
\SPfad(7,7),4\endSPfad
\DuennPunkt(-15,5)
\DuennPunkt(-14,6)
\DuennPunkt(-13,7)
\DuennPunkt(-11,7)
\DuennPunkt(-10,6)
\DuennPunkt(-9,7)
\DuennPunkt(-8,6)
\DuennPunkt(-7,7)
\DuennPunkt(-6,6)
\DuennPunkt(-5,7)
\DuennPunkt(-4,6)
\DuennPunkt(-3,5)
\DuennPunkt(0,5)
\DuennPunkt(1,6)
\DuennPunkt(2,7)
\DuennPunkt(5,7)
\DuennPunkt(6,8)
\DuennPunkt(7,7)
\DuennPunkt(8,6)
\DuennPunkt(9,5)
\Label\o{\longrightarrow}(-1.5,5.5)
\Label\o{\longrightarrow}(-1.5,0.5)
\Label\o{\textrm{{\small $R$}}}(-12,6.8)
\Label\o{\textrm{{\small $F(R)$}}}(3.5,6.8)
\Label\o{\textrm{{\small $Q$ non-primitive;}}}(12,1.5)
\Label\o{\textrm{{\small $Q=\ep$ or ends $DD$}}}(12,0.5)
\Pfad(-14,0),3\endPfad
\Pfad(-11,1),34\endPfad
\Pfad(-7,1),344\endPfad
\Pfad(1,0),3\endPfad
\Pfad(7,1),4\endPfad
\SPfad(-9,1),34\endSPfad
\SPfad(2,1),3\endSPfad
\SPfad(6,2),4\endSPfad
\Label\o{\nearrow}(1.2,6.2)
\Label\o{\swarrow}(1.2,6.2)
\Label\o{\searrow}(7.9,6.2)
\Label\o{\nwarrow}(7.9,6.2)
\Label\o{\textrm{{\small $i$}}}(0.9,6.7)
\Label\o{\textrm{{\small $i$}}}(8.2,6.7)
\Label\l{\longleftarrow}(-8.7,5.4)
\Label\l{\textrm{{\small ---}}}(-7.7,5.4)
\Label\r{\longrightarrow}(-5.4,5.4)
\Label\l{\textrm{{\small ---}}}(-5.1,5.4)
\Label\o{\textrm{{\small $i$}}}(-7,5.0)
\Label\o{\textrm{{\small $Q$}}}(-12,0.8)
\Label\o{\textrm{{\small $F(Q)$}}}(4.5,1.8)
\DuennPunkt(-14,0)
\DuennPunkt(-13,1)
\DuennPunkt(-11,1)
\DuennPunkt(-10,2)
\DuennPunkt(-9,1)
\DuennPunkt(-8,2)
\DuennPunkt(-7,1)
\DuennPunkt(-6,2)
\DuennPunkt(-5,1)
\DuennPunkt(-4,0)
\DuennPunkt(1,0)
\DuennPunkt(2,1)
\DuennPunkt(3,2)
\DuennPunkt(6,2)
\DuennPunkt(7,1)
\DuennPunkt(8,0)
\Label\o{\nearrow}(2.2,1.2)
\Label\o{\swarrow}(2.2,1.2)
\Label\o{\searrow}(6.9,1.2)
\Label\o{\nwarrow}(6.9,1.2)
\Label\o{\textrm{{\small $i$}}}(1.9,1.7)
\Label\o{\textrm{{\small $i$}}}(7.2,1.7)
\Label\l{\longleftarrow}(-9.7,.4)
\Label\l{\textrm{{\small ---}}}(-8.7,.4)
\Label\r{\longrightarrow}(-6.4,.4)
\Label\l{\textrm{{\small ---}}}(-6.1,.4)
\Label\o{\textrm{{\small $i$}}}(-8,0)
\Label\o{ \textrm{\small definition of $F$ on primitive Dyck 
paths}}(0,-2.5)
\]
\vspace*{3mm}

Note that $R=\ep$ in the top left path duplicates a case of the 
bottom left path but no matter: both formulas give the same result.

The map $G$, defined as follows, serves as an inverse of $F$ and 
hence $F$ is indeed a bijection.
Again, $G(\ep)=\ep$ and for a non-primitive Dyck path $P$ with components 
$P_{1},P_{2},\ldots,P_{r}\ (r\ge 2)$, $G(P)=G(P_{1})G(P_{2})\ldots 
G(P_{r})$. By considering the lowest valley vertex, every primitive 
Dyck path has the form $U^{i+1}QD^{i+1}$ with $i\ge 0$ and $Q$ a 
non-primitive Dyck path ($Q=\ep$ in case valley vertices are absent); 
define $G(P)$ by
\[
G(P)=
\begin{cases}
    U UG(R)D (UD)^{i}D & \textrm{if $Q$ ends $UD$, say $Q=RUD$, and} \\
    UG(Q)(UD)^{i}D & \textrm{otherwise.}
\end{cases}
\]

The bijection $F$ is the identity on Dyck paths of size $\le 3$, 
except that it interchanges $U^{3}D^{3}$ and $U^{2}DUD^{2}$. Its 
action on primitive Dyck 4-paths is given in the Figure below.

\Einheit=0.4cm
\[
\Label\o{\longrightarrow}(0,21)
\Label\o{\longrightarrow}(0,16)
\Label\o{\longrightarrow}(0,11)
\Label\o{\longrightarrow}(0,6)
\Label\o{\longrightarrow}(0,1)
\Label\o{\longrightarrow}(0,26)
\Label\o{ \textrm{{\small Dyck path $P$}}}(-5,26)
\Label\o{ \textrm{{\small image $F(P)$}}}(5,26)
\SPfad(-9,0),11111111\endSPfad
\SPfad(1,0),11111111\endSPfad
\SPfad(-9,5),11111111\endSPfad
\SPfad(1,5),11111111\endSPfad
\SPfad(-9,10),11111111\endSPfad
\SPfad(1,10),11111111\endSPfad
\SPfad(-9,15),11111111\endSPfad
\SPfad(1,15),11111111\endSPfad
\SPfad(-9,20),11111111\endSPfad
\SPfad(1,20),11111111\endSPfad
\Pfad(-9,0),33434344\endPfad
\Pfad(1,0),33334444\endPfad
\Pfad(-9,5),33433444\endPfad
\Pfad(1,5),33433444\endPfad
\Pfad(-9,10),33344344\endPfad
\Pfad(1,10),33343444\endPfad
\Pfad(-9,15),33343444\endPfad
\Pfad(1,15),33434344\endPfad
\Pfad(-9,20),33334444\endPfad
\Pfad(1,20),33344344\endPfad
\DuennPunkt(-9,0)
\DuennPunkt(-8,1)
\DuennPunkt(-7,2)
\DuennPunkt(-6,1)
\DuennPunkt(-5,2)
\DuennPunkt(-4,1)
\DuennPunkt(-3,2)
\DuennPunkt(-2,1)
\DuennPunkt(-1,0) 
\DuennPunkt(1,0)
\DuennPunkt(2,1)
\DuennPunkt(3,2)
\DuennPunkt(4,3)
\DuennPunkt(5,4)
\DuennPunkt(6,3)
\DuennPunkt(7,2)
\DuennPunkt(8,1)
\DuennPunkt(9,0)
\DuennPunkt(-9,5)
\DuennPunkt(-8,6)
\DuennPunkt(-7,7)
\DuennPunkt(-6,6)
\DuennPunkt(-5,7)
\DuennPunkt(-4,8)
\DuennPunkt(-3,7)
\DuennPunkt(-2,6)
\DuennPunkt(-1,5) 
\DuennPunkt(1,5)
\DuennPunkt(2,6)
\DuennPunkt(3,7)
\DuennPunkt(4,6)
\DuennPunkt(5,7)
\DuennPunkt(6,8)
\DuennPunkt(7,7)
\DuennPunkt(8,6)
\DuennPunkt(9,5)
\DuennPunkt(-9,10)
\DuennPunkt(-8,11)
\DuennPunkt(-7,12)
\DuennPunkt(-6,13)
\DuennPunkt(-5,12)
\DuennPunkt(-4,11)
\DuennPunkt(-3,12)
\DuennPunkt(-2,11)
\DuennPunkt(-1,10) 
\DuennPunkt(1,10)
\DuennPunkt(2,11)
\DuennPunkt(3,12)
\DuennPunkt(4,13)
\DuennPunkt(5,12)
\DuennPunkt(6,13)
\DuennPunkt(7,12)
\DuennPunkt(8,11)
\DuennPunkt(9,10)
\DuennPunkt(-9,15)
\DuennPunkt(-8,16)
\DuennPunkt(-7,17)
\DuennPunkt(-6,18)
\DuennPunkt(-5,17)
\DuennPunkt(-4,18)
\DuennPunkt(-3,17)
\DuennPunkt(-2,16)
\DuennPunkt(-1,15) 
\DuennPunkt(1,15)
\DuennPunkt(2,16)
\DuennPunkt(3,17)
\DuennPunkt(4,16)
\DuennPunkt(5,17)
\DuennPunkt(6,16)
\DuennPunkt(7,17)
\DuennPunkt(8,16)
\DuennPunkt(9,15)
\DuennPunkt(-9,20)
\DuennPunkt(-8,21)
\DuennPunkt(-7,22)
\DuennPunkt(-6,23)
\DuennPunkt(-5,24)
\DuennPunkt(-4,23)
\DuennPunkt(-3,22)
\DuennPunkt(-2,21)
\DuennPunkt(-1,20) 
\DuennPunkt(1,20)
\DuennPunkt(2,21)
\DuennPunkt(3,22)
\DuennPunkt(4,23)
\DuennPunkt(5,22)
\DuennPunkt(6,21)
\DuennPunkt(7,22)
\DuennPunkt(8,21)
\DuennPunkt(9,20)
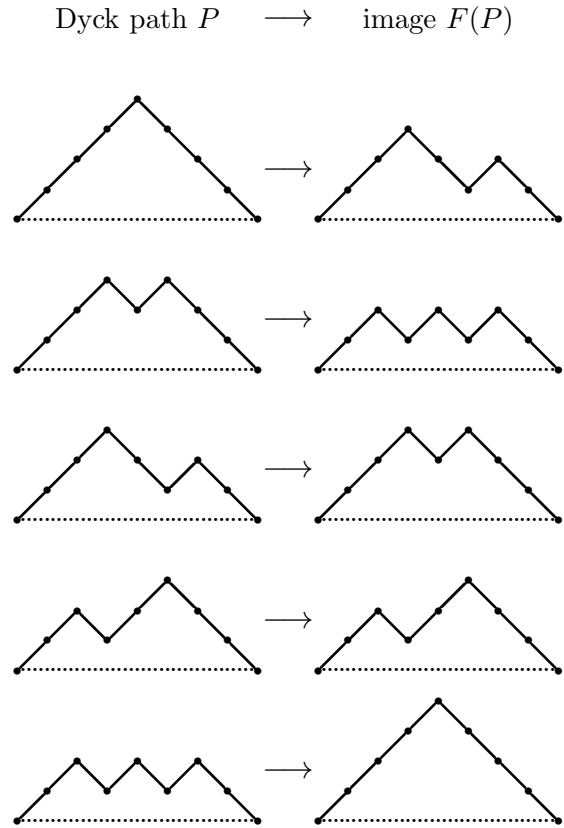
\Label\o{ \textrm{\small action of $F$ on primitive Dyck 4-paths}}(0,-4)
\]

\vspace*{10mm}

{\Large \textbf{4 \ Restriction to \emph{DUU}-avoiding Paths}  }\quad
To analyze the structure of $F$ a key property, clear by induction, 
is that it preserves $\#\,DUU$s, in particular, it 
preserves the property ``path avoids $DUU$''. A 
$DUU$-avoiding Dyck $n$-path corresponds to a composition 
$\c=(c_{1},c_{2},\ldots,c_{h})$ of $n$ via $c_{i}=$ number of $D$s 
ending at height $h-i,\ i=1,2,\ldots,h$ where $h$ is the height of 
the path:

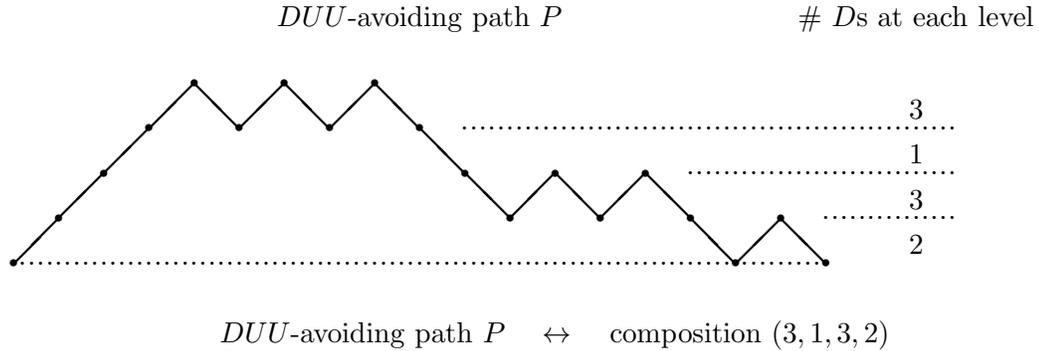
\Einheit=0.6cm
\[
\SPfad(-10,0),111111111111111111\endSPfad
\SPfad(0,3),11111111111\endSPfad
\SPfad(5,2),111111\endSPfad
\SPfad(8,1),111\endSPfad
\Pfad(-10,0),333343434443434434\endPfad
\DuennPunkt(-10,0)
\DuennPunkt(-9,1)
\DuennPunkt(-8,2)
\DuennPunkt(-7,3)
\DuennPunkt(-6,4)
\DuennPunkt(-5,3)
\DuennPunkt(-4,4)
\DuennPunkt(-3,3)
\DuennPunkt(-2,4)
\DuennPunkt(-1,3) 
\DuennPunkt(0,2) 
\DuennPunkt(1,1)
\DuennPunkt(2,2)
\DuennPunkt(3,1)
\DuennPunkt(4,2)
\DuennPunkt(5,1)
\DuennPunkt(6,0)
\DuennPunkt(7,1)
\DuennPunkt(8,0)
\Label\o{ \textrm{\small $DUU$-avoiding path $P$}}(-1,5)
\Label\o{ \textrm{\small \# $D$s at each level}}(10,5)
\Label\o{ \textrm{\small 3}}(10,3)
\Label\o{ \textrm{\small 1}}(10,2)
\Label\o{ \textrm{\small 3}}(10,1)
\Label\o{ \textrm{\small 2}}(10,0)
\Label\o{ \textrm{\small $DUU$-avoiding path $P\quad \leftrightarrow\quad$ 
composition $(3,1,3,2)$}}(2,-2)
\]
\vspace*{2mm}

Under this correspondence, $F$ acts on compositions of $n$\,: $F$ is the 
identity on compositions of length 1, and for $\c=(c_{i})_{i=1}^{r}$ 
with $r\ge 2,\ F(\c)$ is the concatenation of 
$IncrementLast\big(F(c_{1},\ldots,c_{r-2})\big),\,1^{c_{r-1}-1},\,c_{r}$ 
where $IncrementLast$ means ``add 1 to the last entry'' and the 
superscript refers to repetition. In fact, $F$ can be described 
explicitly on compositions of $n$:
\begin{prop}
    For a composition \c of $n$, $F(\c)$ is given by the following 
    algorithm. For each entry $c$ in even position measured from 
    the end $($so the last entry is in position $1)$, replace it by 
    $c-1\ 1$s and increment its left neighbor.
    \label{X}
\end{prop}
For example, 
$4\,2\,1\,5\,2\,3 = \overset{6}{4}\,\overset{5}{2}\,\overset{4}{1}\,\overset{3}{5}\,
\overset{2}{2}\,\overset{1}{3} \rightarrow 1\ 1^{3}\ 3 \ 1^{0}\ 6 \ 
1^{1}\ 3 = 1^{4}\,3\,6\,1\,3$. \qed

Primitive $DUU$-avoiding Dyck $n$-paths correspond to compositions of $n$ that end 
with a 1. Let $\cc_{n}$ denote the set of such compositions. Thus
$\v \cc_{1} \v = 1$ and for $n\ge 2$, \ 
$\v \cc_{n} \v = 2^{n-2}$ since there are $2^{n-2}$ compositions of 
$n-1$.

Denote the length of a composition \c by $\#\c$. The \emph{size} of \c is 
the sum of its entries. The \emph{parity} of \c is the parity (even/odd) 
of $\#\c$. There are two operations on nonempty compositions 
that increment (that is, increase by 1) the size: $P=$ prepend 1, 
and $I=$ increment first entry. For example, for $\c=(4,1,1)$ we 
have size(\c) = 6, $\ \#\c=3,$ the parity of \c is odd, $P(\c)=(1,4,1,1),\ I(\c)=(5,1,1)$. 
\begin{lemma}
    \label{A}
    $P$ changes the parity of a composition while $I$ preserves it.  \qed
\end{lemma}
We'll call $P$ and $I$ \emph{augmentation operators} on $\cc_{n}$ and for 
$A$ an augmentation operator, $A'$ denotes the other one.
\begin{lemma}
    Let $A$ be an augmentation 
    operator. On a composition $\c$ with $\#\c \ge 2$, 
    $A \circ F = F \circ A$ if $\,\#\c$ is odd and $A \circ F = 
    F \circ A'$ if $\,\#\c$ is even.
    \label{B}
\end{lemma}
This follows from Proposition \ref{X}. \qed

Using Lemma \ref{B}, an $F$-orbit $(\c_{1},\ldots,\c_{m})$ in 
$\cc_{n}$ together with an augmentation operator $A_{1} \in\{P,I\}$ 
yields part of an $F$-orbit in $\cc_{n+1}$ via a ``commutative diagram'' 
as shown:
\[
\begin{CD}
    \c_{1} @>F>> \c_{2} @>F>> \ldots @>F>> \c_{i} @>F>> \c_{i+1} @>F>> 
    \ldots @>F>>\c_{m} @>F>> \c_{1} \\
    @VVA_{1}V @VVA_{2}V  @.  @VVA_{i}V @VVA_{i+1}V @. 
    @VVA_{m}V @VVA_{m+1}V \\
    \d_{1} @>F>> \d_{2} @>F>> \ldots @>F>> \d_{i} @>F>> \d_{i+1} @>F>> 
    \ldots @>F>>\d_{m} @>F>> \d_{m+1}
\end{CD}
\]

Let $B(\c_{1},A_{1})$ denote the sequence of compositions $(\d_{1},\ldots,\d_{m})$ 
thus produced. By Lemma \ref{B}, $A_{i+1}= A_{i}$ or $A_{i}'$ 
according as $\,\#\c_{i}$ is odd or even ($1\le i \le m$). Hence, if the orbit of 
$\c_{1}$ contains an even number of compositions of even parity, 
then $A_{m+1}=A_{1}$ and so $\d_{m+1}=\d_{1}$ and $B(\c_{1},A_{1})$ 
is a complete $F$-orbit in $\cc_{n+1}$ for each of $A_{1}=P$ and 
$A_{1}=I$. On the other hand, if the orbit of 
$\c_{1}$ contains an odd number of compositions of even parity, 
then $A_{m+1}=A_{1}'$ and the commutative diagram will extend for 
another $m$ squares before completing an orbit in $\cc_{n+1}$, consisting of the 
concatenation of $B(\c_{1},P)$ and $B(\c_{1},I)$, denoted 
$B(\c_{1},P,I)$. In the former case orbit size is preserved; in the 
latter it is doubled.

Our goal here is to generate $F$-orbits recursively and to get 
induction going, we now need to investigate the parities of the compositions 
comprising these ``bumped-up'' orbits $B(\c,A)$  and 
$B(\c,P,I)$. 
A bit sequence is a sequence of 0s and 1s. \textbf{In the sequel all 
operations on bit sequences are modulo 2}. Let $\s$ denote the partial 
sum operator on bit sequences: $\s\big( (\ep_{1},\ep_{2},\ldots,\ep_{m}) 
\big) =(\ep_{1},\ep_{1}+\ep_{2},\ldots,\ep_{1}+\ep_{2}+\ldots+\ep_{m})$. 
Let $\e_{m}$ denote the all 1s bit sequence of length $m$ and let 
$\e$ denote the infinite sequences of 1s.
Thus $\s\e=(1,0,1,0,1,\ldots)$. 
Let $P$ denote the infinite matrix whose $i$th row ($i\ge 0$) is 
$\s^{i}\e$ ($\s^{i}$ denotes the $i$-fold composition of $\s$). The $(i,j)$ entry $p_{ij}$ of $P$ satisfies $p_{ij}=p_{i-1,j}+p_{i,j-1}$ 
and hence $P$ is the symmetric Pascal matrix mod 2 with $(i,j)$ 
entry =$\:\binom{i+j}{i}$ mod 2. The following lemma will be crucial.
\begin{lemma}
    Fix $k\ge 1$ and let $P_{k}$ denote the $2^{k}\times 2^{k}$ upper 
    left submatrix of $P$. Then the sum modulo $2$ of row $i$ in $P_{k}$ is $0$ 
    for $0\le i \le 2^{k}-1$ and is $1$ for $i=2^{k}-1$.
    \label{P}
\end{lemma}
\textbf{Proof} \quad The sum of row $i$ in $P_{k}$ is, modulo 2,
\[
\sum_{j=0}^{2^{k}-1} p_{ij} = 
\sum_{j=0}^{2^{k}-1}\binom{i+j}{i}=\binom{i+2^{k}}{i+1}=\binom{i+2^{k}}{i+1,2^{k}-1}
\]
and for $i<2^{k}-1$ there is clearly at least one carry in the 
addition of $i+1$ and $2^{k}-1$ in base 2 so that, by Kummer's well 
known criterion, $2\,\v\,\binom{i+2^{k}}{i+1,2^{k}-1}$ and the sum of row $i$
is 0 (mod 2). On the other hand, for $i=2^{k}-1$ there are no 
carries, so $2\nmid \binom{i+2^{k}}{i+1,2^{k}-1}$ and the sum of row $i$ is 1 (mod 2). \qed

Now let $p(\c)$ denote the mod-2 parity of a composition $\c:\ p(\c)=1$ 
if $\,\#\c$ is odd, $=0$ if $\,\#\c$ is even. For purposes of addition 
mod 2, represent the augmentation operators $P$ and $I$ by 0 and 1 
respectively so that, for example, $p(A(\c))=p(\c)+A+1$ for $A=P$ or 
$I$ by Lemma \ref{A}. Then 
the parity of $\d_{i+1}$ above can be obtained from the following 
commutative diagram (all addition modulo 2)
\[
\begin{CD}
   \qquad p(\c_{i})\qquad @>>> p(\c_{i+1}) \\
    @VVAV @VVp(\c_{i})+A+1V  \\
   \qquad \ldots\qquad @>>> p(\c_{i+1})+p(\c_{i})+A
\end{CD}
\]
This leads to
\begin{lemma}
    Let $p_{i}$ denote the parity of $\c_{i}$ so that 
    $\a=(p_{i})_{i=1}^{m}$ is the parity vector for the $F$-orbit 
    $(\c_{i})_{i=1}^{m}$ of the composition $\c_{1}$. Then the parity 
    vector for $B(\c,A)$ is
    \[
    \s\a+\s\e_{m}+(A+1)\e_{m}. 
    \]\qed
\end{lemma}  

Now we are ready to prove the main result of this section concerning 
the orbits of $F$ on primitive $DUU$-avoiding Dyck $n$-paths identified with 
the set $\cc_{n}$ of compositions of $n$ that end with a 1. The parity
of an orbit is the sum mod 2 of the parities of the compositions comprising 
the orbit, in other words, the parity of the total number of entries 
in all the compositions.
\begin{theorem}
    For each $n\ge 1$,
    \begin{itemize}
        \item[$($i\,$)$] all $F$-orbits on $\cc_{n}$ have the same 
        length and this length is a power of $2$.
    
        \item[$($ii\,$)$] all $F$-orbits on $\cc_{n}$ have the same parity.
    
        \item[$($iii\,$)$]  the powers in $($i\,$)$ and the parities 
        in $($ii\,$)$ are given as 
        follows:
	
	For $n=1$, the power $($i.e. the exponent$)$ is $0$ and the parity 
	is $1$.
	For $n=2$, the power and parity are both $0$.
	As $n$ increases from $2$, the powers remain unchanged and the parity 
	stays $0$ except that when $n$ hits a number of the form $2^{k}+1$, the 
	parity becomes $1$, and at the next number, $2^{k}+2$, the power 
	increases by $1$ and the parity reverts to $0$.
    \end{itemize}
\end{theorem}
\textbf{Proof}\quad We consider orbits generated by the augmentation 
operators $P$ and $I$. No orbits are missed because all compositions, 
in particular those  
ending 1, can be generated from the unique composition of 1 by 
successive application of $P$ and $I$. The base cases $n=1,2,3$ are 
clear from the orbits $(1)\to (1),\ (1,1)\to (1,1),\ (2,1) \to 
(1,1,1)\to (2,1)$. To establish the induction step, suppose given an 
orbit, orb($\c$), in $\cc_{2^{k}+1}\ (k\ge 1)$ with parity vector 
$\a=(a_{i})_{i=1}^{2^{k}}$ and (total) parity 1. Then the next orbit 
$B(\c,P,I)$ has parity vector
\[
\u_{1}=(\s\, \a,\s\, \a+\e_{2^{k}})+\s\,\e_{2^{k+1}}
\]
with parity ($\s\,\a$'s cancel out) $\underbrace{1+1+\ldots+1}_{2^{k}}+\underbrace{1+0+1+0+\ldots+1+0}_{2^{k+1}}=0$ for $k\ge 1$.
Successively ``bump up'' this orbit using $A=\ep_{1},\ep_{2},\ldots,$ 
in turn until the parity hits 1 again. 
With Sum$(\vv)$ denoting the sum 
of the entries in \vv, the successive parity 
vectors $\u_{1},\u_{2},\ldots$ are given by 
\begin{multline*}
\u_{i}=\big(\s^{i}\a,\s^{i}\a+\sum_{j=1}^{i-2}\textrm{Sum}(\s^{j}\a)\s^{i-1-j}\e_{2^{k}} + \s^{i-1}\e_{2^{k}}\big) + \\
\s^{i}\e_{2^{k+1}} + \s^{i-1}\e_{2^{k+1}} 
+\sum_{j=1}^{i-2}\ep_{j}\s^{i-1-j}\e_{2^{k+1}} + (\ep_{i-1}+1)\e_{2^{k+1}}.
\end{multline*}

Applying Lemma \ref{P} we see that, independent of the $\ep_{i}$'s, $\u_{i}$ has sum 0 for 
$i<2^{k}-1$ and sum 1 for $i=2^{k}-1$. This establishes the 
induction step in the theorem. \qed

\begin{cor}
    For $n\ge 2$, the length of each $F$-orbit in $\p_{n}(DUU)$ is $2^{k}$ where $k$ 
    is the number of bits in the base-$2$ expansion of $n-2$.
    \label{base2}
\end{cor}
\textbf{Proof}\quad This is just a restatement of part of the 
preceding Theorem. \qed

\vspace*{10mm}

{\Large \textbf{5 \ The Orbits of $\mathbf{F}$}  }\quad 
The preceding section analyzed $F$ on $\p(DUU)$, paths avoiding 
$DUU$. Now we consider $F$ on $\p[DUU]$, the primitive Dyck paths containing a 
$DUU$. Every $P \in \p[DUU]$ has the form $AQB$ where 
\begin{itemize}
    \item[(i)] $A$ consists of one or more $U$s

    \item[(ii)] $C:=AB \in \p(DUU)$

    \item[(iii)] $Q \notin \p $ and $Q$ ends $DD$ (and hence $Q$ contains 
    a $DUU$ at \gl).
\end{itemize}

To see this, locate the rightmost of the lowest $DUU$s in $P$, say at 
height $h$. Then $A=U^{h},\ Q$ starts at step number $h+1$ and 
extends through the matching downstep of the middle $U$ in this 
rightmost lowest $DUU$, and $B$ consists of the rest of the path.

\Einheit=0.4cm
\[
\Pfad(-15,0),3\endPfad
\Pfad(-13,2),33\endPfad
\Pfad(-8,4),433\endPfad
\Pfad(-2,5),44\endPfad
\Pfad(6,3),4\endPfad
\Pfad(14,1),4\endPfad
\SPfad(-14,1),3\endSPfad
\SPfad(-11,4),413\endSPfad
\SPfad(-5,5),413\endSPfad
\SPfad(13,2),4\endSPfad
\DuennPunkt(-15,0)
\DuennPunkt(-14,1)
\DuennPunkt(-13,2)
\DuennPunkt(-12,3)
\DuennPunkt(-11,4) 
\DuennPunkt(-8,4)
\DuennPunkt(-7,3)
\DuennPunkt(-6,4)
\DuennPunkt(-5,5)
\DuennPunkt(-2,5)
\DuennPunkt(-1,4)
\DuennPunkt(0,3)
\DuennPunkt(6,3)
\DuennPunkt(7,2)
\DuennPunkt(13,2)
\DuennPunkt(14,1)
\DuennPunkt(15,0)
\red{
\DuennPunkt(1,4)
\DuennPunkt(2,3)
\DuennPunkt(4,3)
\DuennPunkt(5,4)
\DuennPunkt(8,3)
\DuennPunkt(9,2)
\DuennPunkt(11,2)
\DuennPunkt(12,3)
\Pfad(0,3),34\endPfad
\Pfad(4,3),34\endPfad
\Pfad(7,2),34\endPfad
\Pfad(11,2),34\endPfad
\SPfad(2,3),34\endSPfad
\SPfad(9,2),34\endSPfad
}
\blue{
\Pfad(-15,0),111111111111111111111111111111\endPfad
\Pfad(-15,0),2222222\endPfad
\Pfad(-12,0),2222222\endPfad
\Pfad(0,0),2222222\endPfad
\Pfad(15,0),2222222\endPfad
}
\Label\o{\uparrow}(-7,1.7)
\Label\o{\textrm{{\footnotesize $h$}}}(-7,0.8)
\Label\u{\downarrow}(-7,1.2)
\Label\o{\textrm{{\footnotesize $A$}}}(-13.5,6)
\Label\o{\textrm{{\footnotesize $\leftarrow$ matching 
$\rightarrow$}}}(-3.5,2.7)
\Label\o{\textrm{{\footnotesize $Q$}}}(-6,6)
\Label\o{\textrm{{\footnotesize $B$}}}(7,6)
\Label\o{\textrm{{\footnotesize red $UD$s may be absent}}}(7,4.5)
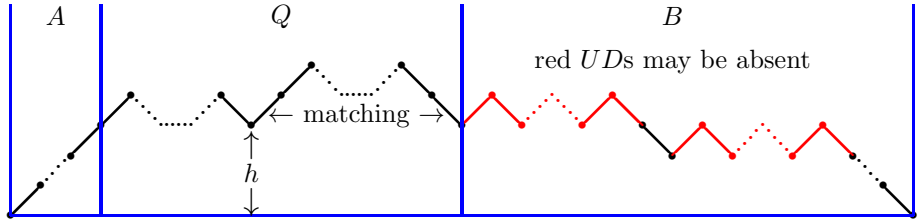
\Label\o{\textrm{{\small The $AQB$ decomposition of a path containing a $DUU$}}}(0,-3)
\]
\vspace*{2mm}

Call the path $AB$ the ($DUU$-avoiding) \emph{skeleton} of $P$ and 
$Q$ the ($DUU$-containing) \emph{body} of $P$. In case $P\in\p(DUU)$, its skeleton 
is itself and its body is empty. If the skeleton of $P$ is $UD$, 
then $P$ is uniquely determined by its skeleton and body. On the 
other hand, a 
skeleton of size $\ge 2$ and a nonempty body determine precisely two 
paths $P$ in $\p[DUU]$, obtained by inserting the body at either the 
top or the bottom of the first peak upstep in the skeleton, as 
illustrated.

\Einheit=0.4cm
\[
\Pfad(-17,0),3334344344\endPfad
\Pfad(-5,0),33\endPfad
\Pfad(-2,2),34344344\endPfad
\Pfad(7,0),333\endPfad
\Pfad(11,3),4344344\endPfad
\SPfad(-17,0),1111111111\endSPfad
\SPfad(-5,0),11111111111\endSPfad
\SPfad(7,0),11111111111\endSPfad
\DuennPunkt(-17,0)
\DuennPunkt(-16,1)
\DuennPunkt(-15,2)
\DuennPunkt(-14,3)
\DuennPunkt(-13,2)
\DuennPunkt(-12,3)
\DuennPunkt(-11,2)
\DuennPunkt(-10,1)
\DuennPunkt(-9,2)
\DuennPunkt(-8,1)
\DuennPunkt(-7,0)
\DuennPunkt(-5,0)
\DuennPunkt(-4,1)
\DuennPunkt(-3,2)
\DuennPunkt(-2,2)
\DuennPunkt(-1,3)
\DuennPunkt(0,2)
\DuennPunkt(1,3)
\DuennPunkt(2,2)
\DuennPunkt(3,1)
\DuennPunkt(4,2)
\DuennPunkt(5,1)
\DuennPunkt(6,0)
\DuennPunkt(7,0)
\DuennPunkt(8,1)
\DuennPunkt(9,2)
\DuennPunkt(10,3)
\DuennPunkt(11,3)
\DuennPunkt(12,2)
\DuennPunkt(13,3)
\DuennPunkt(14,2)
\DuennPunkt(15,1)
\DuennPunkt(16,2)
\DuennPunkt(17,1)
\DuennPunkt(18,0)
\Label\o{\textrm{{\footnotesize $S$}}}(-12,-2)
\Label\o{\textrm{{\footnotesize two possible $P$s}}}(6.5,-2)
\Label\o{\textrm{{\footnotesize $B$}}}(-2.5,1.8)
\Label\o{\textrm{{\footnotesize $B$}}}(10.5,2.8)
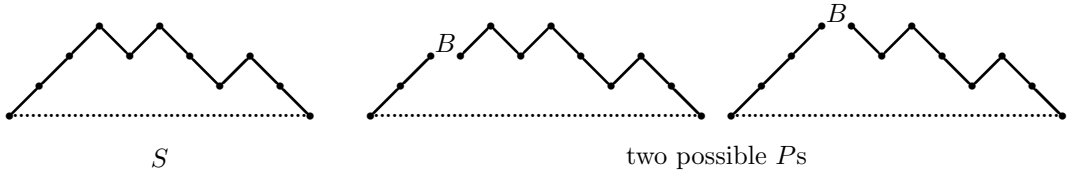
\Label\o{\textrm{{\small Recapturing a path $P\in \p[DUU]$ from a skeleton 
$S$ and body $B$}}}(0,-4)
\]
\vspace*{2mm}

Thus paths in $\p[DUU]$ correspond bijectively to triples $(S,B,pos)$ 
where $S\in\p(DUU)$ is the skeleton, $B\ne \ep$ is the body, and $pos 
=top$ or $bot$ according as $B$ is positioned at the top or bottom 
of the first peak upstep in $S$, with the proviso that $pos=top$ if 
$S=UD$.

In these terms, $F$ can be specified on $\p[DUU]$ as follows.
\begin{prop}
\[
F\big( (S,B,pos)\big)=
\begin{cases}
    (F(S),F(B),\:pos\,) \textrm{ if height$(S)$ is odd, and} \\
    (F(S),F(B),\:pos'\,) \textrm{ if height$(S)$ is even.}
\end{cases}
\]    
\end{prop}
\textbf{Proof}\quad Let $h(P)$ denote the height of the terminal point 
of the lowest $DUU$ in $P\in \p[DUU]$. The result clearly holds for 
$h(P)=1$. If $h(P)\ge 2$, then $P$ has the form 
$U^{2}Q(UD)^{a}D(UD)^{b}D$ with $a,b\ge 0$ and $Q$ a Dyck path that 
ends $DD$. So $F(P)=U^{b+1}F(Q)(UD)^{a+1}D^{b+1}$ and $h(Q)=h(P)-2$. 
These two facts are the basis for a proof by induction that begins as 
follows. If $h(Q)=0$, then the body of $F(P)$ has position = bottom, 
while the body of $P$ has position bottom or top according as $a\ge 
1$ or $a=0$. In the former case, the skeleton of $P$ has height 3 and 
position has been preserved, in the latter height 2 and position has 
been reversed. \qed

Iterating the skeleton-body-position decomposition on each component, 
a Dyck path has a forest representation as illustrated below. Each vertex 
represents a skeleton and is labeled with the corresponding composition. 
When needed, a color ($top$ or $bot$) is also applied to a vertex to 
capture the position of that skeleton's body.

\Einheit=0.4cm
\[
\Pfad(-16,0),33344334443433343433444334344344\endPfad
\SPfad(-16,0),11111111111111111111111111111111\endSPfad
\DuennPunkt(-16,0)
\DuennPunkt(-15,1)
\DuennPunkt(-14,2)
\DuennPunkt(-13,3)
\DuennPunkt(-12,2)
\DuennPunkt(-11,1)
\DuennPunkt(-10,2)
\DuennPunkt(-9,3)
\DuennPunkt(-8,2)
\DuennPunkt(-7,1)
\DuennPunkt(-6,0)
\DuennPunkt(-5,1)
\DuennPunkt(-4,0)
\DuennPunkt(-3,1)
\DuennPunkt(-2,2)
\DuennPunkt(-1,3)
\DuennPunkt(0,2)
\DuennPunkt(1,3)
\DuennPunkt(2,2)
\DuennPunkt(3,3)
\DuennPunkt(4,4)
\DuennPunkt(5,3)
\DuennPunkt(6,2)
\DuennPunkt(7,1)
\DuennPunkt(8,2)
\DuennPunkt(9,3)
\DuennPunkt(10,2)
\DuennPunkt(11,3)
\DuennPunkt(12,2)
\DuennPunkt(13,1)
\DuennPunkt(14,2)
\DuennPunkt(15,1)
\DuennPunkt(16,0)
\]
\begin{center}

\begin{pspicture}(-6,-1.4)(6,3)

\psline(-4,1)(-3,0)(-2,1)
\psline(1,2)(2,1)(3,2)
\psline(2,2)(2,1)(3,0)(4,1)
\rput(-3,0){$\bullet$}
\rput(0,0){$\bullet$}
\rput(3,0){$\bullet$}

\psdots(-4,1)(-2,1)(2,1)(4,1)(1,2)(2,2)(3,2)

\rput(-4.2,1.2){\textrm{{\footnotesize 11}}}
\rput(-1.8,1.2){\textrm{{\footnotesize 11}}}
\rput(2.5,1){\textrm{{\footnotesize 1}}}

\rput(0.9,2.2){\textrm{{\footnotesize 1}}}
\rput(2,2.3){\textrm{{\footnotesize 1}}}
\rput(3.1,2.2){\textrm{{\footnotesize 11}}}

\rput(4.1,1.2){\textrm{{\footnotesize 21}}}

\rput(0,-.3){\textrm{{\footnotesize 1}}}
\rput(-3,-.3){\textrm{{\footnotesize 1}}}
\rput(3,-.3){\textrm{{\footnotesize \quad 11, bot}}}

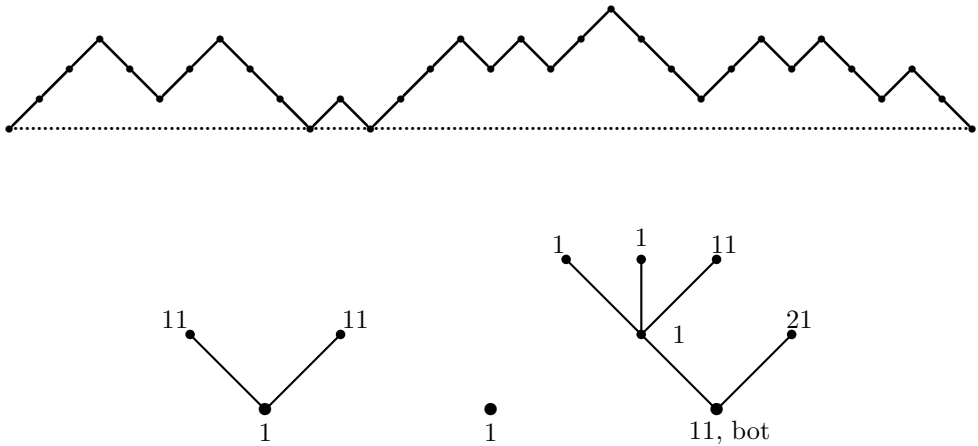
\rput(0,-1.3){\textrm{{\small A Dyck path and corresponding LCO forest}}}

\end{pspicture}
\end{center} 
The 3 trees in the forest correspond to the 3 components of the Dyck 
path. The skeleton of the first component is $UD$ and its body has 2 
identical components, each consisting of a skeleton alone, yielding 
the leftmost tree. The skeleton of the third component is $UUDD$ and 
its body is positioned at the bottom of its first peak upstep, and so 
on.
Call this forest the LCO (labeled, colored, ordered) forest
corresponding to the Dyck path. Here is the precise definition.
\begin{defn}
    An LCO forest is a labeled, colored, ordered forest such that
    \begin{itemize}
    \item  the underlying forest consists of a list of ordered trees 
    (a tree may consist of a root only)

    \item  no vertex has outdegree $1$ $($i.e., exactly one child\,$)$

    \item  each vertex is labeled with a composition that ends $1$

    \item  each vertex possessing children and labeled with a composition 
    of size $\ge 2$ is also colored $top$ or $bot$

    \item  For each leaf $($i.e. vertex with a parent but no child\,$)$ that 
    is the rightmost child of its parent, its label composition has 
    size $\ge 2$.
\end{itemize}
\end{defn}

The \emph{size} of an LCO forest is the sum  of the sizes of its label compositions. 
The correspondence Dyck path $\leftrightarrow$ LCO forest 
preserves size, and primitive Dyck paths correspond to one-tree forests. 
Thus we have 

\begin{prop}
    The number of LCO forests of size $n$ is the Catalan number 
    $C_{n}$, as is the number of one-tree LCO forests of size $n+1$. \qed
\end{prop} 

The $C_{4}=14$ one-tree LCO forests corresponding to primitive 
Dyck 5-paths are shown, partitioned into $F$-orbits. 
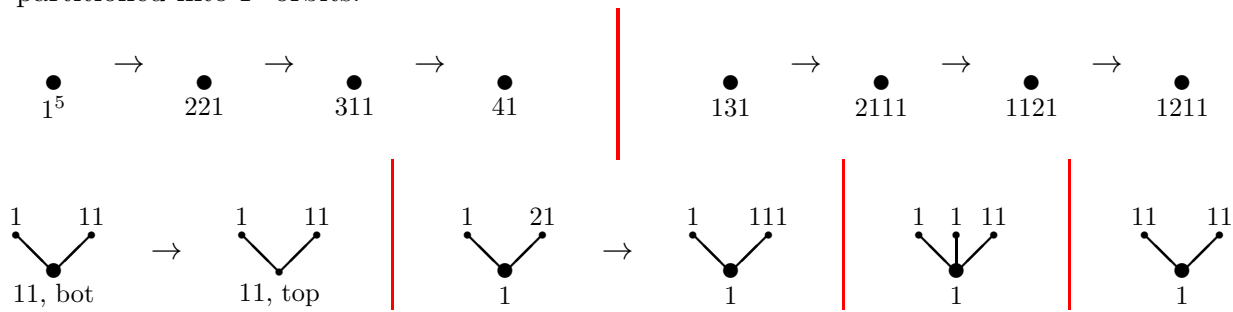
\Einheit=0.5cm
\[
\Label\o{\rightarrow}(-13,5)
\Label\o{\rightarrow}(-9,5)
\Label\o{\rightarrow}(-5,5)
\Label\u{ \textrm{{\footnotesize $1^{5}$}}}(-15,5)
\Label\u{\textrm{{\footnotesize 221}}}(-11,5)
\Label\u{\textrm{{\footnotesize 311}}}(-7,5)
\Label\u{\textrm{{\footnotesize 41}}}(-3,5)
\NormalPunkt(-15,5)
\NormalPunkt(-11,5)
\NormalPunkt(-7,5)
\NormalPunkt(-3,5)
\Label\o{\rightarrow}(13,5)
\Label\o{\rightarrow}(9,5)
\Label\o{\rightarrow}(5,5)
\Label\u{\textrm{{\footnotesize 1211}} }(15,5)
\Label\u{\textrm{{\footnotesize 1121}} }(11,5)
\Label\u{\textrm{{\footnotesize 2111}} }(7,5)
\Label\u{\textrm{{\footnotesize 131}} }(3,5)
\NormalPunkt(15,5)
\NormalPunkt(11,5)
\NormalPunkt(7,5)
\NormalPunkt(3,5)
\Label\o{\rightarrow}(-12,0)
\Label\o{\rightarrow}(-0,0)
\Label\u{ \textrm{{\footnotesize 11,\ bot}}}(-15,0)
\Label\u{ \textrm{{\footnotesize 11,\ top}}}(-9,0)
\Label\o{\textrm{{\footnotesize 1}} }(-16,1)
\Label\o{\textrm{{\footnotesize 11}} }(-14,1)
\Label\o{\textrm{{\footnotesize 1}} }(-10,1)
\Label\o{\textrm{{\footnotesize 11}} }(-8,1)
\Label\o{\textrm{{\footnotesize 1}} }(-4,1)
\Label\o{\textrm{{\footnotesize 21}} }(-2,1)
\Label\o{\textrm{{\footnotesize 1}} }(2,1)
\Label\o{\textrm{{\footnotesize 111}} }(4,1)
\Label\o{\textrm{{\footnotesize 1}} }(8,1)
\Label\o{\textrm{{\footnotesize 1}} }(9,1)
\Label\o{\textrm{{\footnotesize 11}} }(10,1)
\Label\o{\textrm{{\footnotesize 11}} }(16,1)
\Label\o{\textrm{{\footnotesize 11}} }(14,1)
\Label\u{\textrm{{\footnotesize 1}}}(-3,0)
\Label\u{\textrm{{\footnotesize 1}}}(3,0)
\Label\u{\textrm{{\footnotesize 1}}}(9,0)
\Label\u{\textrm{{\footnotesize 1}}}(15,0)
\red{\Pfad(0,3),2222\endPfad
\Pfad(-6,-1),2222\endPfad
\Pfad(6,-1),2222\endPfad
\Pfad(12,-1),2222\endPfad}
\Pfad(-16,1),4\endPfad
\Pfad(-15,0),3\endPfad
\Pfad(-10,1),4\endPfad
\Pfad(-9,0),3\endPfad
\Pfad(-4,1),4\endPfad
\Pfad(-3,0),3\endPfad
\Pfad(2,1),4\endPfad
\Pfad(3,0),3\endPfad
\Pfad(8,1),4\endPfad
\Pfad(9,0),3\endPfad
\Pfad(9,0),2\endPfad
\Pfad(14,1),4\endPfad
\Pfad(15,0),3\endPfad
\DuennPunkt(-16,1)
\NormalPunkt(-15,0)
\DuennPunkt(-14,1)
\DuennPunkt(-10,1)
\DuennPunkt(-9,0)
\DuennPunkt(-8,1)
\DuennPunkt(-4,1)
\NormalPunkt(-3,0)
\DuennPunkt(-2,1)
\DuennPunkt(2,1)
\NormalPunkt(3,0)
\DuennPunkt(4,1)
\DuennPunkt(8,1)
\DuennPunkt(9,1)
\NormalPunkt(9,0)
\DuennPunkt(10,1)
\DuennPunkt(14,1)
\DuennPunkt(16,1)
\NormalPunkt(15,0)
\Label\u{\textrm{{\small The LCO one-tree forests of size 5, partitioned into 
$F$-orbits}} }(0,-2)
\]
\vspace*{1mm}

We can now give an explicit description of $F$ on Dyck paths 
identified with LCO forests. On an LCO forest, $F$ acts as follows:
\begin{itemize}
    \item  the underlying list of ordered trees is preserved

    \item  each label $\c$ becomes $F(\c)$ as defined in Prop.\:\ref{X}

    \item  each color ($top$/$bot$) is preserved or switched according 
    as the associated label \c has odd or even length.
\end{itemize}

From this description and Cor. \ref{base2}, the size of the $F$-orbit of a Dyck path $P$ 
can be determined as follows. In the LCO forest for $P$, let $\ell$ 
denote the maximum size of a leaf label and $i$ the maximum size of an 
internal (i.e., non-leaf) label (note that an isolated root is an 
internal vertex). Let $k$ denote the number of bits in the base-2 
expansion of $\max\{\ell-2,i-1\}$. Then the $F$-orbit of $P$ has size 
$2^{k}$.

It is also possible to specify orbit sizes in terms of subpath 
avoidance. For Dyck paths $Q$ and $R$, let $Q$ \emph{top} $R$ (resp. $Q$ \emph{bot} 
$R$) denote the Dyck path obtained by inserting $R$ at the top (resp. 
bottom) of the first peak upstep in $Q$. Then the $F$-orbit of a Dyck path $P$ has 
size $\le 2^{k}$ iff $P$ avoids subpaths in the set $\{Q\ top\ R,\ Q \ 
bot\ R\, :\, R\ne \ep,\ Q\in\p_{i}(DUU),\ 2^{k-1}+1 < i \le 
2^{k}+1\}$. For $k\ge 1,$ listing these $Q$s explicitly would give 
$2^{2^{k}}-2^{2^{k-1}}$ proscribed patterns of the form $Q\ top\ R,\ 
R\ne \ep$ (and the same number of the form $Q\ bot\ R$). For $k=0$, 
that is, for fixed points of $F$, the proscribed patterns are 
$UP^{+}UDD$ and $UUP^{+}DD$ with $P^{+}$ a nonempty Dyck path,
and avoiding the first of these amounts to avoiding the subpath $DUDD$.

The \gf\ for the number of $F$-orbits of size $\le 2^{k}$ can be found 
using the ``symbolic'' method \cite{flaj}. With 
$F_{k}(x),\: G_{k}(x),\: H_{k}(x)$ denoting the respective \gf s for 
general Dyck paths, primitive Dyck paths, and primitive Dyck paths 
that end $DD$ ($x$ always marking size), we find
\begin{eqnarray*}
    F_{k}(x) & = & 1+G_{k}(x)F_{k}(x)  \\
    G_{k}(x) & = & x 
    +\frac{x(1-(2x)^{2^{k}}}{1-2x}\big(x+(F_{k}(x)-1)H_{k}(x)\big)  \\
    H_{k}(x) & = & G_{k}(x)-x
\end{eqnarray*}
leading to 
\[
F_{k}(x)=\frac{1-a_{k}-\sqrt{1-4x-\frac{\textrm{{\small $a_{k}(2-a_{k})x $}}}{\textrm{{\small $1-x$}}}}}{2x-a_{k}},
\]
where $a_{k}=(2x)^{2^{k}+1}$. In this formulation it is clear, as 
expected, that $\lim_{k \to 
\infty}F_{k}(x)=\frac{1-\sqrt{1-4x}}{2x}$, the \gf\ for the Catalan 
numbers. The counting sequence for fixed points of $F$, with \gf\ 
$F_{0}(x)$, is sequence 
\htmladdnormallink{A086625}{http://www.research.att.com:80/cgi-bin/access.cgi/as/njas/sequences/eisA.cgi?Anum=A086625}
in 
\htmladdnormallink{OEIS}{http://www.research.att.com/~njas/sequences/Seis.html} .

\vspace*{10mm}

{\Large \textbf{6 \ An Application}  }\quad 
Ordered trees and binary trees are manifestations of the Catalan 
numbers
\htmladdnormallink{A000108}{http://www.research.att.com:80/cgi-bin/access.cgi/as/njas/sequences/eisA.cgi?Anum=A000108} . 
Donaghey \cite{motz77,restricted77} lists several types of restricted tree 
counted by the Motzkin numbers
\htmladdnormallink{A001006}{http://www.research.att.com:80/cgi-bin/access.cgi/as/njas/sequences/eisA.cgi?Anum=A001006} . 
In particular, the following result 
is implicit in item III\,C of \cite{restricted77}.
\begin{prop}
    The Motzkin number $M_{n}$ counts right-planted binary trees on 
    $n+1$ edges with no erasable vertices.
    \label{erasable}
\end{prop}
Here, planted means the root has only one child, and erasable 
refers to a vertex incident with precisely 2 edges \emph{both of the same 
slope}---the vertex could then be erased, preserving the slope, to produce a smaller binary 
tree. The $M_{3}=4$ such trees on 4 edges are shown.

\Einheit=0.6cm
\[
\Pfad(-7,2),3\endPfad
\Pfad(-7,2),43\endPfad
\Pfad(-7,0),3\endPfad
\Pfad(-3,0),33\endPfad
\Pfad(-3,2),4\endPfad
\Pfad(-2,3),4\endPfad
\Pfad(1,3),43\endPfad
\Pfad(2,2),4\endPfad
\Pfad(2,0),3\endPfad
\Pfad(5,0),3\endPfad
\Pfad(5,2),4\endPfad
\Pfad(5,2),3\endPfad
\Pfad(5,4),4\endPfad
\DuennPunkt(-7,0) 
\DuennPunkt(-7,2)
\DuennPunkt(-6,1)
\DuennPunkt(-6,3)
\DuennPunkt(-5,2)
\DuennPunkt(-3,0)
\DuennPunkt(-3,2)
\DuennPunkt(-2,1)
\DuennPunkt(-2,3)
\DuennPunkt(-1,2) 
\DuennPunkt(1,3)
\DuennPunkt(2,0)
\DuennPunkt(2,2)
\DuennPunkt(3,1)
\DuennPunkt(3,3)
\DuennPunkt(5,0)
\DuennPunkt(5,4)
\DuennPunkt(5,2)
\DuennPunkt(6,1)
\DuennPunkt(6,3)
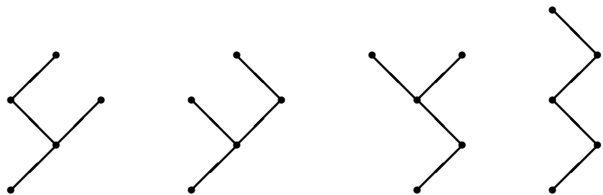
\Label\o{ \textrm{\small The right-planted binary 4-trees with no erasable 
vertices}}(0,-2)
\]

\vspace*{2mm}

Translated to Dyck paths, Prop. \ref{erasable} is equivalent to
\begin{prop}
    $M_{n}$ counts Dyck $(n+1)$-paths that end $DD$ and avoid 
    subpaths $DUDU$ and $UUP^{+}DD$ with $P^{+}$ denoting a nonempty Dyck subpath.
    \label{UUXDD}
\end{prop}
We will use $F$ to give a bijective proof of Prop. \ref{UUXDD} based 
on the fact \cite{udu} that $M_{n}$ also counts $DUD$-avoiding Dyck 
$(n+1)$-paths. (Of course, path reversal shows that $\#\,UDU$s and 
$\#\,DUD$s are equidistributed on Dyck paths.) Define statistics $X$ and $Y$ on Dyck paths by 
$X=\#\:DUD$s and $Y=\#\:DUDU$s $ +\ \#\:UUP^{+}DD$s + [paths ends with 
$UD$] (Iverson notation) so that the paths in Prop. \ref{UUXDD} are 
those with $Y=0$. Prop. \ref{UUXDD} then follows from
\begin{prop}
    On Dyck $n$-paths with $n\ge 2$, $F$ sends the statistic $X$ to 
    the statistic $Y$.
\end{prop}
\textbf{Proof}\quad Routine by induction from the recursive definition 
of $F$. However, using the explicit form of $F$, it is also possible to 
specify precisely which $DUD$s correspond to each of the three summands 
in $Y$. For this purpose, given a $DUD$ in a Dyck path $P$, say $D_{1}U_{2}D_{3}$ 
(subscripts used simply to  
identify the individual steps), let 
$\ss(D_{1}U_{2}D_{3})$ 
denote the longest Dyck subpath of $P$ containing $D_{1}U_{2}D_{3}$ in its 
skeleton and let $h$ denote the height at which $D_{1}U_{2}D_{3}$ terminates 
in $\ss(D_{1}U_{2}D_{3})$. If $h$ is odd, $D_{1}U_{2}D_{3}$ is 
immediately followed in $P$ by $D_{4}$ or by $UD_{4}$ (it cannot be followed by 
$UU$). In either case, let $U_{4}$ be the matching upstep for 
$D_{4}$. Then the steps $D_{1},U_{2},D_{3},U_{4}$ show up in $F(P)$ as 
part of 
a subpath $U_{4}U_{2}P^{+}D_{3}D_{4}$ with $P^{+}$ a Dyck path that 
ends $D_{1}$. On the other hand, if $h$ is even, $D_{1}U_{2}D_{3}$ 
either (i) ends the path (here $\ss(D_{1}U_{2}D_{3})=P$ and $h=0$) or is 
immediately followed by (ii) $U_{4}$ or (iii) $D$. In case (iii), 
let $U_{4}$ be the matching upstep. Then $D_{1},U_{2},D_{3},U_{4}$ show 
up in $F(P)$ as a subpath in that order (cases (ii) and (iii)) or 
$F(P)$ ends $U_{2}D_{3}$ (case (i)). The details are left to the reader.

\vspace*{10mm}

{\Large \textbf{7 \ Statistics Suggested by LCO Forests}  }\quad 
There are various natural statistics on LCO forests, some of which 
give interesting counting results. Here we present two such. First
let us count one-tree LCO forests by size of root label. This is equivalent 
to counting primitive Dyck paths by skeleton size. Recall that 
the generalized Catalan number sequence $\big(C^{(j)}_{n}\big)_{n\ge 
0}$ with $C^{(j)}_{n}:=\frac{j}{2n+j}\binom{2n+j}{n}$ is the $j$-fold 
convolution of the ordinary Catalan number sequence
\htmladdnormallink{A000108}{http://www.research.att.com:80/cgi-bin/access.cgi/as/njas/sequences/eisA.cgi?Anum=A000108}.
(See \cite{woan} for a nice bijective proof.)
And, as noted above, in the skeleton-body-position decomposition of a 
primitive Dyck path, if the body is nonempty it contains a $DUU$ at 
(its own)
\gl \ and ends $DD$.
\begin{lemma}
    The number of Dyck $n$-paths that  contain a $DUU$ at 
\gl \ and end $DD$ is $C^{(4)}_{n-3}$.
\end{lemma}
\textbf{Proof}\quad In such a path, let $U_{0}$ denote the middle $U$ 
of the \emph{last} $DUU$ at \gl. The path then has the form $AU_{0}BD$ 
where $A$ and $B$ are arbitrary \emph{nonempty} Dyck paths, counted 
by $C^{(2)}_{n-1}$. So the desired counting sequence is the convolution 
of $\big(C^{(2)}_{n-1}\big)$ with itself and, taking the $U_{0}D$ into 
account, the lemma follows. \qed

The number of primitive $DUU$-avoiding Dyck $k$-paths is 1 if $k=1$, 
and $2^{k-2}$ if $k\ge 2$. But if $k\ge 2$, there are two choices 
(top/bottom) to insert the body. So the number of primitive Dyck $(n+1)$-paths 
with skeleton size $k$ is $2^{k-1}C^{(4)}_{n-k-2}$ for $1\le k \le 
n-2$ and is $2^{n-1}$ for $k=n+1$. Since there are $C_{n}$ primitive Dyck 
$(n+1)$-paths altogether, we have established the following  identity.
\begin{prop}
   \[
C_{n} = 2^{n-1} + \sum_{k=1}^{n-2}\frac{2^{k}}{n-k}\binom{2n-2k}{n-2-k}.
\] 
\end{prop} \qed

Lastly, turn an LCO forest into an LCO tree by joining all roots to a new root. 
The purpose of doing this is so that isolated roots in the forest will 
qualify as leaves in the tree. The symbolic method then yields 
\begin{prop}
  The \gf\ for LCO trees by number of leaves $(x$ marks size, $y$ 
  marks number of leaves\,$)$ is
  \[
  \frac{1-\sqrt{1-4x\:\frac{\textrm{{\small $1-x$}}}{\textrm{{\small $1-xy$}}}}}{2x}.
  \]
\end{prop}
The first few values are given in the following table.

\[
\begin{array}{c|cccccccc}
	n^{\textstyle{\,\backslash \,k}} &  1 & 2 & 3 & 4 & 5 & 6 & 7 & 8 \\
\hline 
	1&    1 &   & & &  & & & \\
 	2&    1 & 1 & & &  & & &   \\
	3&    2 & 2  & 1 & &  & & &  \\
	4&    4 & 6  & 3 & 1 &  & & &   \\ 
	5&    8 & 17 & 12 & 4 & 1 & & &   \\
	6&    16 & 46 & 44 & 20 & 5 & 1 & &   \\
	7&    32 & 120 & 150 & 90 & 30  & 6 & 1 &   \\
        8&    64 & 304 & 482 & 370 & 160  & 42 & 7 & 1   \\ 
 \end{array}
\]
\centerline{{\small number of LCO trees of size $n$ with $k$ leaves}}


\begin{thebibliography}{99}
    
    
    \bibitem{twobij04} David Callan, Two bijections for Dyck path parameters,
\htmladdnormallink{math.CO/0406381}{http://front.math.ucdavis.edu/math.CO/0406381}, 2004, 4pp.

\bibitem{catfine} David Callan, Some bijections and identities for the Catalan and Fine numbers,
\htmladdnormallink{S\'{e}m. Lothar. Combin.}{http://www.math.ethz.ch/EMIS/journals/SLC/index.html}
\textbf{53} (2004/06), Art. B53e, 16 pp. 

\bibitem{invol1999} Emeric Deutsch, An involution on Dyck paths and its consequences. \emph{Discrete Math.} \textbf{204} (1999), no. 1-3, 163--166. 
 
\bibitem{bij1998} Emeric Deutsch, A bijection on Dyck paths and its consequences,
\emph{Discrete Math.} \textbf{179} (1998), no. 1-3, 253--256. 

\bibitem{ordered} Emeric Deutsch, A bijection on ordered trees 
and its consequences, \emph{J. Combin. Theory Ser. A} \textbf{90} (2000), no. 1, 210--215. 

\bibitem{simple2003} Emeric Deutsch and Sergi Elizalde, A simple and unusual bijection for Dyck paths and its consequences,
\emph{Ann. Comb.} \textbf{7} (2003), no. 3, 281--297.  

\bibitem{don80} Robert Donaghey,  Automorphisms on Catalan trees and bracketings. \emph{J. Combinatorial Theory Ser. B}
\textbf{29} (1980), no. 1, 75--90. MR0584162 

\bibitem{acp44}
Donald Knuth, \emph{Art of Computer Programming, Vol.4, Fascicle 4: Generating all Trees -- History 
of Combinatorial Generation}, Addison-Wesley, 2006, vi+120pp, draft available from 
\htmladdnormallink{http://www-cs-faculty.stanford.edu/$\,\widetilde{\ }\,$knuth/fasc4a.ps.gz}{http://www-cs-faculty.stanford.edu/~knuth/fasc4a.ps.gz}   
    
\bibitem{lalanne92} J.-C. Lalanne,  Une involution sur les chemins de 
Dyck, \emph{European J. Combin.} \textbf{13} (1992), no. 6, 477--487.

\bibitem{lalanne93} J.-C. Lalanne, Sur une involution sur les chemins de Dyck, Conference on Formal Power Series and Algebraic 
Combinatorics \emph{Theoret. Comput. Sci.} \textbf{117} (1993), no. 1-2, 203--215. 

\bibitem{vaille97} J. Vaill\'{e}, Une bijection explicative de plusieurs 
propri\'{e}t\'{e}s remarquables des ponts, 
\emph{European J. Combin.} \textbf{18} (1997), no. 1, 117--124. 



    



\bibitem{motz77} Robert Donaghey and Louis Shapiro, Motzkin numbers, 
\emph{ J. Combinatorial Theory Ser. A}  \textbf{23}, 291--301, 1977. MR0505544 

\bibitem{restricted77} Robert Donaghey, Restricted plane tree representations 
of four Motzkin-Catalan equations,  \emph{J. Combinatorial Theory 
Ser. B} \textbf{22}, (1977), no. 2, 114--121, 1977.  MR0432532


\bibitem{udu} Y. Sun, The statistic ``number of udu's'' in Dyck paths, \emph{Disc. Math.},
              \textbf{287} (2004), Issue 1-3 (October 2004), 177-186.

\bibitem{flaj} Robert Sedgewick and Philippe Flajolet, 
An Introduction to the Analysis of Algorithms,
Addison-Wesley, 1996.

\bibitem{woan} Wen-jin Woan,
Uniform partitions of lattice paths and Chung-Feller generalizations. 
\emph{Amer. Math. Monthly} 108 (2001), no. 6, 556--559.

\end{thebibliography}
\end{document}